\documentclass[12pt,english]{article}
\usepackage[T1]{fontenc}
\usepackage[latin9]{inputenc}
\usepackage{color}
\usepackage{mathrsfs}
\usepackage{amsmath}
\usepackage{amssymb}

\makeatletter

\usepackage[active]{srcltx}

\usepackage{tikz}
\usepackage{tikz-cd}
\usetikzlibrary{shapes.geometric}
\usepackage{ bbold }

\usepackage{a4wide}
\usepackage[english]{babel}
\usepackage{amsthm}
\usepackage{amsopn}
\usepackage[active]{srcltx}
\usepackage[notref,notcite]{showkeys}
\usepackage{hyperref}
\usepackage[T1]{fontenc}

\newtheorem{theorem}{Theorem}

\usepackage[notref,notcite]{showkeys}

\usepackage{hyperref}

\usepackage{lastpage}
\usepackage{fancyhdr}
\makeatother

\usepackage{babel}
\begin{document}
\title{Special properties of Gossez's example}
\author{M.D. Voisei}
\date{{}}
\maketitle
\begin{abstract}
In this note new properties for the Gossez example are presented in
regard to its representability, closedness, and maximal monotonicity
with respect to the two dual systems it naturally inhabits. 
\end{abstract}

\strut


\section{Introduction}

In \cite[p. 89]{MR0298492} Gossez presents his famous example of
a monotone operator which exhibits strange range properties (see also
\cite{MR0397485,Bauschke2020}), namely, $G:\ell_{1}\to\ell_{\infty}=\ell_{1}^{*}$
given by 
\begin{equation}
(Gx)_{n}:=-\sum_{k<n}x_{k}+\sum_{k>n}x_{k}=\sum_{k\in\mathbb{N}}x_{k}\alpha(k,n),\ x=(x_{k})_{k\in\mathbb{N}}\in\ell_{1},\ n\in\mathbb{N},\label{eq:}
\end{equation}
where $\alpha(k,n)=\left\{ \begin{array}{lll}
-1, & \mathrm{if} & k<n\\
0, & \mathrm{if} & k=n\\
1, & \mathrm{if} & k>n
\end{array}\right.$, $k,n\in\mathbb{N}$. 

Here
\[
\ell_{1}:=\{x=(x_{n})_{n\in\mathbb{N}}\mid\|x\|_{\ell_{1}}:=\sum_{n\in\mathbb{N}}|x_{n}|<+\infty\},\ \ell_{\infty}:=\{x=(x_{n})_{n\in\mathbb{N}}\mid\|x\|_{\ell_{\infty}}:=\sup_{n\in\mathbb{N}}|x_{n}|<+\infty\}
\]
\[
\mathfrak{c}:=\{y=(y_{n})_{n\in\mathbb{N}}\mid\lim_{n\to\infty}y_{n}\in\mathbb{R}\}.
\]

The linear operator $G$ is well-defined and bounded since, for every
$x=(x_{n})_{n}\in\ell_{1}$, 
\[
\|Gx\|_{\ell_{\infty}}=\sup_{n\in\mathbb{N}}|(Gx)_{n}|\le\sum_{n\in\mathbb{N}}|x_{n}|=\|x\|_{\ell_{1}}.
\]

For a topologically vector space $(E,\mu)$, $A\subset E$, and $f,g:E\rightarrow\overline{\mathbb{R}}$
we denote by $E^{*}$ the topological dual of $E$ and by $[f\le g]$
$:=\{x\in E\mid f(x)\leq g(x)\}$; the sets $[f=g]$, $[f<g]$, $[f>g]$,
$[f\ge g]$ are similarly defined. 

A \emph{dual system} is a triple $(X,Y,c)$ consisting of two vector
spaces $X$, $Y$ and a bilinear map $c=\langle\cdot,\cdot\rangle:X\times Y\to\mathbb{R}$.
The weakest topology on $X$ for which all linear forms $\{X\ni x\to\langle x,y\rangle\mid y\in Y\}$
are continuous is denoted by $\sigma(X,Y)$. This topology is locally
convex and called \emph{the weak topology of} $X$ with respect to
the duality $(X,Y,c)$. Similarly one defines the weak topology $\sigma(Y,X)$
on $Y$. 

For $S\subset X$, 
\[
S^{\perp}:=\{y\in Y\mid\forall x\in S,\ \langle x,y\rangle=0\}
\]
denotes the \emph{orthogonal} \emph{of} $S$. Similarly, for $A\subset Y$,
$A^{\perp}:=\{x\in X\mid\forall y\in A,\ \langle x,y\rangle=0\}$
is the \emph{orthogonal} of $A$. 

A template for dual systems $(X,Y,\langle\cdot,\cdot\rangle)$ is
given by $(X,X^{*},\langle\cdot,\cdot\rangle)$ where $(X,\tau)$
is a locally convex space, $X^{*}$ is the topological dual of $(X,\tau)$,
and $\langle x,x^{*}\rangle=x^{*}(x)$, $x\in X$, $x^{*}\in X^{*}$. 

In this paper we study several operators and sets defined in special
dual systems, e.g., $(\ell_{1},\ell_{\infty})$, $(\ell_{\infty},\ell_{\infty}^{*})$,
$(\ell_{1}\times\ell_{\infty},\ell_{\infty}\times\ell_{1})$, or $(\ell_{1}\times\ell_{\infty},\ell_{\infty}\times\ell_{\infty}^{*})$;
$\ell_{1}^{*}$ is identified with $\ell_{\infty}$, and recall that
the coupling for the duality $(\ell_{1},\ell_{\infty})$ is
\[
\langle x,y\rangle_{\ell_{1}\times\ell_{\infty}}:=\langle x,y\rangle:=\sum_{n\in\mathbb{N}}x_{n}y_{n},\ x=(x_{n})_{n\in\mathbb{N}}\in\ell_{1},y=(x_{n})_{n\in\mathbb{N}}\in\ell_{\infty},
\]
and it coincides with the coupling of $(\ell_{\infty}^{*},\ell_{\infty})$,
that is, 
\[
\forall x\in\ell_{1},\ \forall y\in\ell_{\infty},\ \langle x,y\rangle_{\ell_{\infty}^{*}\times\ell_{\infty}}=\langle x,y\rangle_{\ell_{1}\times\ell_{\infty}}
\]
which allows us to use the simplified notation of $\langle\cdot,\cdot\rangle$
for both couplings. 

Given a dual system $(X,Y,c=\langle\cdot,\cdot\rangle)$, $Z:=X\times Y$
forms naturally a dual system $(Z,Z,\cdot)$ where, for $z=(x,y),\ w=(u,v)\in Z$,
\begin{equation}
z\cdot w:=c(x,v)+c(u,y).\label{eq:-161}
\end{equation}
The space $Z$ is endowed with a locally convex topology $\mu$ compatible
with the natural duality $(Z,Z,\cdot)$, that is, for which $(Z,\mu)^{*}=Z$;
for example $\mu=\sigma(Z,Z)=\sigma(X,Y)\times\sigma(Y,X)$. With
respect to the natural dual system $(Z,Z,\cdot)$, the conjugate of
$f$ is denoted by 
\begin{equation}
f^{\square}:Z\rightarrow\overline{\mathbb{R}},\quad f^{\square}(z)=\sup\{z\cdot z^{\prime}-f(z^{\prime})\mid z^{\prime}\in Z\}.\label{eq:-160}
\end{equation}
By the biconjugate formula, $f^{\square\square}=\operatorname*{cl}\operatorname*{conv}f$
whenever $f^{\square}$ (or $\operatorname*{cl}\operatorname*{conv}f$)
is proper. 

To $A\subset Z$ we associate $c_{A}:Z\rightarrow\overline{\mathbb{R}}$,
$c_{A}:=c+\iota_{A}$, where $\iota_{A}(z)=0$, for $z\in A$, $\iota_{A}(z)=+\infty$,
otherwise; the \emph{Fitzpatrick function} of $A$ 
\begin{equation}
\varphi_{A}:Z\rightarrow\overline{\mathbb{R}},\ \varphi_{A}(z):=c_{A}^{\square}(z)=\sup\{z\cdot w-c(w)\mid w\in A\},\label{eq:-9}
\end{equation}
and $\psi_{A}:Z\rightarrow\overline{\mathbb{R}}$, $\psi_{A}:=\varphi_{A}^{\square}=c_{A}^{\square\square}$. 

In particular when $\varphi_{A}$ is proper (or $\operatorname*{cl}\operatorname*{conv}c_{A}$
is proper), 
\[
\psi_{A}=\operatorname*{cl}\operatorname*{conv}c_{A},\ \varphi_{A}=\psi_{A}^{\square}.
\]
Note that $\varphi_{A}$ is improper iff $\varphi_{A}$ is identically
$-\infty$ (when $A=\emptyset$) or $\varphi_{A}$ is identically
$+\infty$. 

Similarly, for a multi\-function $T:X\rightrightarrows Y$ we define
$c_{T}:=c_{\operatorname*{Graph}T}$, $\psi_{T}:=\psi_{\operatorname*{Graph}T}$,
and the \emph{Fitzpatrick function} of $T$, $\varphi_{T}:=\varphi_{\operatorname*{Graph}T}$.
By convention $\varphi_{\emptyset}=-\infty$, $c_{\emptyset}=\operatorname*{conv}c_{\emptyset}=\psi_{\emptyset}=+\infty$
in agreement with the usual conventions of $\inf_{\emptyset}=+\infty$,
$\sup_{\emptyset}=-\infty$. 

In expanded form, for $T:X\rightrightarrows Y$ and $z=(x,y)\in Z$,
$(x,y)\in X\times Y$

\vspace{-.5cm}

\begin{align}
\varphi_{T}(z) & =\sup\{z\cdot w-c(w)\mid w\in\operatorname*{Graph}T\},\nonumber \\
\varphi_{T}(x,y) & =\sup\{\langle x,v\rangle+\langle u,y\rangle-\langle u,v\rangle\mid(u,v)\in\operatorname*{Graph}T\}.\label{def-Ff}
\end{align}

In the sequel, for $A\subset X\times Y$ we denote by 
\begin{equation}
\neg A:=\{(x,y)\in X\times Y\mid(x,-y)\in A\}.\label{eq:-2}
\end{equation}

For a multi-valued operator $T:X\rightrightarrows Y$ we define $\neg T:X\rightrightarrows Y$
via $\operatorname*{Graph}(\neg T):=\neg\operatorname*{Graph}(T)$,
or, equivalently, 
\[
(\neg T)(x)=-T(x),\ x\in D(\neg T):=D(T).
\]

Consider the following classes of functions on $Z=X\times Y$:

\medskip

\noindent $\mathscr{C}:=\mathscr{C}(Z):=\{f\in\Lambda(Z)\mid f\geq c\}$
the class of proper convex functions that are greater than the coupling
$c$ on $Z$; 

\medskip

\noindent $\mathscr{R}:=\mathscr{R}(Z):=\Gamma(Z)\cap\mathscr{C}(Z)$
the class of proper convex lower semicontinuous functions that are
greater than the coupling $c$ on $Z$; 

\medskip

\noindent $\mathscr{D}:=\mathscr{D}(Z):=\{f\in\mathscr{R}(Z)\mid f^{\square}\geq c\}$
the class of proper convex lower semicontinuous functions $f$ such
that $f$, $f^{\square}$ are greater than the coupling $c$ on $Z$;\textcolor{red}{{} }

\smallskip

A multi-function $T:X\rightrightarrows Y$ is 
\begin{itemize}
\item \emph{monotone} if, for all $y_{1}\in T(x_{1})$, $y_{2}\in T(x_{2})$,
$\left\langle x_{1}-x_{2},y_{1}-y_{2}\right\rangle \geq0$; 
\item \emph{maximal monotone} if $\operatorname*{Graph}T$ is maximal in
the sense of graph inclusion as monotone subsets of $X\times Y$; 
\item \emph{unique} in $Z=X\times Y$ if $T$ is monotone and admits a unique
maximal monotone extension in $Z$. 
\item \emph{representable} in $Z=X\times Y$ if $\operatorname*{Graph}T=[f=c]$
for some $f\in\mathscr{R}(Z)$; in this case $f$ is called a \emph{representative}
of $T$. We denote by $\mathscr{R}_{T}:=\mathscr{R}_{T}(Z)$ the class
of representatives of $T$; 
\item \emph{dual--representable} in $Z$ if $\operatorname*{Graph}T=[f=c]$
for some $f\in\mathscr{D}(Z)$; in this case $f$ is called a \emph{d--representative}
of $T$. We denote by $\mathscr{D}_{T}:=\mathscr{D}_{T}(Z)$ the class
of d--representatives of $T$; 
\item of \emph{negative infimum type }(\emph{NI} for short)\emph{ in $Z$}
if $\varphi_{T}\geq c$ (in $Z$);
\end{itemize}
A subset $A\subset Z$ is said to have a certain property if $A=\operatorname*{Graph}T$
for some muti-valued $T:X\rightrightarrows Y$ which has that same
property. 

Gossez's example in (\ref{eq:}) distinguishes these notions in the
sense seen in Theorem \ref{sds} below. 

Recall that $T:X\rightrightarrows Y$ is maximal monotone iff it is
representable and NI iff it is dual-representable and unique (see
e.g. \cite{MR2207807,MR2389004,MR2453098,MR2577332,MR2594359,MR3917361}). 

\section{Preliminaries}

\begin{theorem} \label{G} \emph{(i)} $R(G)\subset\mathfrak{c}$
and $G$ is one-to-one. 

\medskip

\emph{(ii)} With respect to the duality $(\ell_{1},\ell_{\infty})$
the operator $G$ is linear bounded skew, that is, for every $x\in D(G)=\ell_{1}$,
$\langle x,Gx\rangle=0$. Moreover
\begin{equation}
\forall x,y\in\ell_{1},\ \langle x,Gy\rangle=-\langle y,Gx\rangle,\label{eq:-8}
\end{equation}
or, equivalently, $(\operatorname*{Graph}G)^{\perp}=\operatorname*{Graph}G$
(relative to the duality $(\ell_{1}\times\ell_{\infty},\ell_{\infty}\times\ell_{1})$).
\end{theorem} 

Consider $G^{*}:\ell_{\infty}^{*}\to\ell_{\infty}$ the adjoint of
$G:\ell_{1}\to\ell_{\infty}$, given by 
\[
\langle y,G^{*}\mu\rangle=\langle\mu,Gy\rangle,\ y\in\ell_{1},\ \mu\in\ell_{\infty}^{*}.
\]

In order to describe $G^{*}$, in this article $\ell_{\infty}$ is
identified with $C(\beta\mathbb{N})$ the space of continuous functions
over the Hausdorff separated compact space $\beta\mathbb{N}$ which
is the Stone-\v{C}ech compactification of $\mathbb{N}$ in the sense
that every $x\in\ell_{\infty}$ is identified with and uniquely extended
to $\beta x\in C(\beta\mathbb{N})$ such that $\beta x|_{\mathbb{N}}=x$;
$\ell_{\infty}^{*}$ is identified with $C(\beta\mathbb{N})^{*}$
which, by the Riesz's Representation Theorem, is identified with $M(\beta\mathbb{N})$
the space of signed regular Borel measures of finite total variation
on $\beta\mathbb{N}$ in the sense that for every $F\in C(\beta\mathbb{N})^{*}$
there is a unique $\mu\in M(\beta\mathbb{N})$ such that 
\[
\forall f\in C(\beta\mathbb{N}),\ F(f)=\int_{\beta\mathbb{N}}fd\mu,
\]
$F$ is denoted by $\mu$ and, for $\mu\in\ell_{\infty}^{*}$, $x\in\ell_{\infty}$,
we write $\mu(x)=:\langle\mu,x\rangle={\displaystyle \int_{\beta\mathbb{N}}}\beta xd\mu$. 

\begin{theorem} \label{G*} \emph{(i)} $G^{*}\mu=-\mu(\beta\mathbb{N}\smallsetminus\mathbb{N})\mathbb{1}-G\bar{\mu},\ \mu\in M(\beta\mathbb{N})$,
where $\bar{\mu}:=(\mu(\{k\}))_{k\in\mathbb{N}}\in\ell_{1}$,
\begin{equation}
\operatorname*{Ker}(G^{*}):=(G^{*})^{-1}(0)=\{\mu\in M(\beta\mathbb{N})\mid\mu(\beta\mathbb{N}\smallsetminus\mathbb{N})=0,\ \forall n\in\mathbb{N},\ \mu(\{n\})=0\},\label{eq:-7}
\end{equation}
and $G^{*}$ is not one-to-one. 

\medskip

\emph{(ii)} $R(G)$ is $\sigma(\ell_{\infty},\ell_{1})-$dense and
neither (strongly-)closed nor dense in $\ell_{\infty}$. Here $\overline{R(G)}$
denotes the strong closure of $R(G)$. \end{theorem}

\section*{Dual-system results}

\begin{theorem} \label{fds} With respect to the duality $(\ell_{1},\ell_{\infty})$
the operator $G:\ell_{1}\to\ell_{\infty}$ is maximal monotone with
\begin{equation}
\varphi_{G}=\psi_{G}=\iota_{\operatorname*{Graph}G}=c_{G}.\label{eq:-10}
\end{equation}
Similar considerations hold for $\neg G$. \end{theorem}

\begin{theorem} \label{sds} \emph{(i)} With respect to the duality
$(\ell_{\infty}^{*},\ell_{\infty})$, that is, when we see 
\[
G:D(G)=\ell_{1}\subsetneq\ell_{1}^{**}=\ell_{\infty}^{*}\to\ell_{\infty},
\]
$G$ is skew monotone unique and not NI, not representable, not $\sigma(\ell_{\infty}^{*},\ell_{\infty})\times\sigma(\ell_{\infty},\ell_{\infty}^{*})-$closed,
not dual-representable, and not maximal monotone. The only maximal
monotone extension of $G$ is $\neg G^{*}$ and (relative to the duality
$(\ell_{1}\times\ell_{\infty},\ell_{\infty}\times\ell_{\infty}^{*})$)
\[
\operatorname*{Graph}G\subset(\operatorname*{Graph}G)^{\perp}=\operatorname*{Graph}(\neg G^{*}).
\]

The Fitzpatrick function of $G$ with respect to the duality $(\ell_{\infty}^{*},\ell_{\infty})$
is 
\begin{equation}
\Phi_{G}=\iota_{\operatorname*{Graph}(\neg G^{*})},\ \Psi_{G}:=\Phi_{G}^{\square}=\iota_{\operatorname*{cl}(\operatorname*{Graph}(G))},\label{eq:-11}
\end{equation}
and 
\[
L:=\operatorname*{cl}(\operatorname*{Graph}G)=[\Phi_{G}=c]=[\Psi_{G}=c]
\]
is a skew (that is, $L\subset[c=0]$) linear subspace of $\ell_{\infty}^{*}\times\ell_{\infty}$
which is unique representable but not NI, not dual-rep\-re\-sent\-a\-ble,
and not maximal monotone; but with $\Phi_{G}=\Phi_{L}$, $\Psi_{G}=\Psi_{L}$,
and 
\begin{equation}
\operatorname*{Graph}G\subsetneq L=[\Phi_{L}=c]=[\Psi_{L}=c]\subsetneq\operatorname*{Graph}(G^{*}).\label{eq:-27}
\end{equation}

\emph{(ii)} In contrast the operator $\neg G:\ell_{1}\subsetneq\ell_{\infty}^{*}\to\ell_{\infty}$,
$(\neg G)(x)=-Gx$, $x\in D(\neg G)=\ell_{1}$, is NI but not maximal
monotone with unique maximal monotone extension the skew operator
with graph $\operatorname*{cl}(\operatorname*{Graph}(\neg G))=\neg L$. 

Here ``$\operatorname*{cl}$'' stand for the closure with respect
to any topology compatible with the duality $(\ell_{\infty}^{*}\times\ell_{\infty},\ell_{\infty}\times\ell_{\infty}^{*})$
such as $\sigma(\ell_{\infty}^{*},\ell_{\infty})\times\sigma(\ell_{\infty},\ell_{\infty}^{*})$.
\end{theorem}

\bibliographystyle{plain}

\end{document}